# Efficiency of the maximum partial likelihood estimator for nested case control sampling

LARRY GOLDSTEIN[1] and HAIMENG ZHANG[2]

[1]*Kaprielian Hall, Room 108, 3620 Vermont Avenue, Los Angeles, California, 90089-2532, USA.
E-mail: larry@math.usc.edu*

[2]*Department of Mathematics and Statistics, Mississippi State University, Mississippi State, MS 39762. E-mail: HZhang@math.msstate.edu*

In making inference on the relation between failure and exposure histories in the Cox semiparametric model, the maximum partial likelihood estimator (MPLE) of the finite dimensional odds parameter, and the Breslow estimator of the baseline survival function, are known to achieve full efficiency when data is available for all time on all cohort members, even when the covariates are time dependent. When cohort sizes become too large for the collection of complete data, sampling schemes such as nested case control sampling must be used and, under various models, there exist estimators based on the same information as the MPLE having smaller asymptotic variance.

Though the MPLE is therefore not efficient under sampling in general, it approaches efficiency in highly stratified situations, or instances where the covariate values are increasingly less dependent upon the past, when the covariate distribution, not depending on the real parameter of interest, is unknown and there is no censoring. In particular, in such situations, when using the nested case control sampling design, both the MPLE and the Breslow estimator of the baseline survival function achieve the information lower bound both in the distributional and the minimax senses in the limit as the number of cohort members tends to infinity.

## 1. Introduction

For many epidemiologic studies, the cohort from which failures are observed is simply too large for the collection of full exposure data, and in order to make inference on the connection between exposure history and failure it becomes a matter of practical necessity to sample. For a cohort followed over time, one of the simplest sampling schemes, termed nested case control sampling [15], is to choose a fixed number of controls to compare to the failure at each failure time. Though it has previously been shown that the maximum partial likelihood estimator (MPLE) in the Cox semi-parametric model achieves full efficiency when data is available for all time on all cohort members, the same is no longer true in certain situations when schemes such as nested case control sampling are used. In counterpoint to such cases, here we explore a model where the MPLE is efficient, in both the distributional and minimax senses, for the nested case control sampling scheme. We also show that similar remarks apply as well to the Breslow estimator of the baseline hazard. Knowing in which situations the MPLE is close to efficient provides some guidelines on when it may be applied with little risk of efficiency loss, and when other estimators, perhaps depending on additional modeling assumptions, should be considered as an alternative.

In the standard Cox model [5], a common but unspecified baseline hazard function $\lambda(t)$ is assumed to apply to all cohort members. The relation between exposure and failure is the one of most interest, and is modeled by the real parameter $\theta$ specifying the increased relative risk, having the exponential form $e^{\theta Z}$, say, for an individual with covariate $Z$. The unknown baseline is considered for the most part to be a nuisance parameter. When covariate information is available on all cohort members, the maximum partial likelihood estimator (MPLE) makes inference on the parametric component of such models by maximizing a 'partial likelihood', that is, the product of the conditional probabilities, over all failures $i_j$, that individual $i_j$ failed given that the individuals $\mathcal{R}_{i_j}$ were also at risk to fail when $i_j$ failed,

$$L(\theta) = \prod_{i_j} \frac{e^{\theta Z_{i_j}}}{\sum_{k \in \mathcal{R}_{i_j}} e^{\theta Z_k}}. \tag{1}$$

We note that the unspecified baseline hazard cancels upon forming this conditional probability.

When data is only available on some sampled subset $\widetilde{\mathcal{R}}_{i_j}$ of the entire cohort $\mathcal{R}_{i_j}$, an estimator may be formed by replacing $\mathcal{R}_{i_j}$ by $\widetilde{\mathcal{R}}_{i_j}$, (see [4]), possibly then mandating the use of weights so that the MPLE remains consistent. Nested case control sampling, which does not require the use of such weights, is the instance where $\widetilde{\mathcal{R}}_{i_j}$ consists of the failure $i_j$ and $m-1$ non-failed individuals to serve as controls, chosen uniformly at random for those at risk at the time of the failure.

One price to pay for the ability to estimate $\theta$ while leaving the nonparametric baseline hazard unspecified, and the subsequent use of the MPLE, is that it is not a true likelihood being maximized, and efficiency concerns arise. In particular, it is not clear whether one can construct estimators that depend on the same data as the MPLE but have better performance. In the paper of Begun *et al.* [1], however, these concerns are put to rest in the full cohort case where the covariates are time fixed, as the authors demonstrate that in that situation the MPLE achieves the semi-parametric efficiency bound. Greenwood and Wefelmeyer [10] show the MPLE is efficient in the full cohort situation even when the covariates are allowed to depend on time. Similar remarks also apply to the Breslow estimator of the baseline hazard.

The situation is different under sampling: Robins *et al.* [14] has shown that for time fixed covariates the MPLE is not efficient under nested case control sampling. In this situation, there may exist modified estimators that take advantage of the time fixed nature of the covariates, in that the exposure for a control sampled in the past is still valid at a future failure time. In time varying covariate models, Chen [6] among others, have modified the MPLE to yield consistent estimators of the parametric parameter that have smaller asymptotic variance than the MPLE. The estimator proposed in [6] uses covariates sampled for other failures at time points near to that of a given failure to take advantage of already available information. Here, to realize a practical efficiency benefit, the sequence of failure times must be sufficiently dense and the covariates must not be varying too rapidly in time. Though in the time fixed covariate situation the modified estimator uses information from the past specifically, in both cases one relies on the dependence of the covariate values over time to realize some efficiency gain; for the time varying covariate models, such modified estimators will perform better the stronger the time dependence. Due to the various improvements on the performance of the MPLE, it becomes less clear in just which ways its performance can be improved, or, in other words, whether the MPLE fails to be efficient for reasons in addition to the ones by which these modified estimators achieve their gains.

Showing that there is some sense in which the MPLE for nested case control sampling is efficient is therefore valuable for two reasons. First, it limits the scope of the search for estimators that might improve the MPLE's performance. Second, it indicates the use of the simple MPLE, and not a more complex version of same, in situations that achieve or approximate those in which it cannot be improved.

Based on the known instances where the MPLE fails to be efficient under sampling, to find models where it is, by contrast, efficient, we are led to consider situations where covariate information collected for one failure is not useful at any other failure time. Indeed, such situations are fairly common in epidemiologic studies, in particular, when highly stratified cohorts are followed over a short period of time. Due to the short time under study, the covariates may be considered time fixed, and there is, for that same reason, little or no censoring. Last, in such cases, the groups corresponding to the terms in the product of the partial likelihood are independent, or very nearly so. A continuous time covariate model where the failures are spaced far apart relative to the correlation time of the covariates will also have the property that the covariate values at one failure time will be nearly independent of those at any other. In fact, in the limit, this latter situation becomes the former, highly stratified case. Thus we are led to a time fixed covariate model $f$ having no censoring, where we observe $n$ independent units of information, each consisting of the observed failure from a cohort of a possibly random number $\eta$ of individuals who are comparable to the failure, the covariate value of the failure, and the covariate values of $m-1$ sampled controls.

A concrete example of such a situation is the study of occupational exposure to electromagnetic fields, or EMF and leukemia [11], which is fairly typical of cancer registry based case-control studies. The cohort is the adult male population in mid-Sweden followed over 1983–1987 for cases of leukemia. Two controls were sampled from risk sets based on the age of the 250 leukemia cases, matching on year of birth and geographic location. In this study, with the four-year follow-up and fine stratification, there is little censoring and almost all strata have at most one failure, thus the sampling model considered here very closely approximates the circumstances of the study.

It is easy to verify that in these situations, letting $Z$ be the distribution of the i.i.d. covariates, under the null $\theta_0 = 0$ the information $-E[\partial^2 \log L(\theta)/\partial \theta^2] = \sigma_{\text{MPLE}}^{-2}$, where

$$\sigma_{\text{MPLE}}^{-2} = \left(\frac{m-1}{m}\right) \text{Var}(Z),$$

where $L(\theta)$ is as in (1) with the set of those at risk $\mathcal{R}_{i_j}$ replaced by the nested case control sampled risk set $\widetilde{\mathcal{R}}_{i_j}$. Hence, under regularity (see, e.g., [4, 7, 8]) the MPLE $\hat{\theta}_n$ is asymptotically normal and satisfies

$$\sqrt{n}(\hat{\theta}_n - \theta_0) \to \mathcal{N}(0, \sigma_{\text{MPLE}}^2).$$

Our main result, Theorem 2.5, shows that when considering a growing cohort size, the limiting effective information in the data, $I_*(\theta_0)$, equals $\sigma_{\text{MPLE}}^{-2}$, and that the MPLE is efficient in the limit in both the convolution lower bound and minimax senses. Theorem 2.6 shows similar remarks apply to the Breslow estimator of the baseline survival function.

When the complete set of covariate values is observed it is unimportant whether the covariate distribution is considered known or unknown. Again, the situation when sampling is different; knowing the covariate distribution allows one to estimate large sample quantities with some accuracy. Consequently, the hypothesis of Theorem 2.5 includes the assumption that the covariate distribution is unknown, and the subsequent analysis must therefore handle two infinite dimensional nuisance parameters, one for the unknown baseline density, the other for the unknown covariate distribution. In particular, the results leave open the possibility of improved estimators that take advantage of a known covariate distribution. Nevertheless, such improvements must necessarily depend on having information about, and correctly modeling, the covariate distribution, and consequently invite the possibility of bias due to modeling misspecification.

We consider the Cox model under the usual exponential relative risk, though the methods here may be applied for other relative risk forms, as was accomplished in [10] for the full cohort, time varying covariate model. The methods here also extend to accommodate censoring, though this generalization requires the inclusion of a third infinite dimensional parameter, the censoring density and consequently the handling of an additional operator corresponding to the unknown censoring density.

The outline of this work is as follows: In Section 2.1 we review and slightly modify the theory in [1] for the calculation of information bounds in semi-parametric models to accommodate a pair of unknown densities. In Section 2.2 we further specialize that theory to the case at hand and formally state our model and the main results that were outlined above. Application of the theory presented in Section 2.1 for the relative risk parameter $\theta$ requires verification of three assumptions. The first, Assumption 2.1, is that certain collections of perturbations form a subspace. The second, Assumption 2.2, is connected to the Hellinger differentiability of the observation density $f$, in particular, that perturbations of the nonparametric baseline and covariate density affect $f$ by amounts given by operators $A$ and $B$ evaluated on the respective perturbations, and that perturbing the parametric parameter results in a score $\rho_0$. The third, Assumption 2.3, is that the orthogonal projection of the parametric score $\rho_0$ is contained in a certain subspace, $\mathbb{K}$. In order to proceed as quickly as possible to the calculation of the information bounds in Section 4, we present in Section 3 only a subset of the properties eventually required of the operators $A$ and $B$ and of the score $\rho_0$.

The remaining properties required of $A$ and $B$ are shown in the Appendix in Sections A.1 and A.2. An outline of the verification of Assumptions 2.1, 2.2 and 2.3 is given in Section A.3; the detailed calculations can be found in the technical report [9]. Remarks on the modifications made to the theory in [1] that are necessary for our application can be found in Section A.4.

## 2. Information bounds for sampling in the Cox model

In Section 2.1 we review and adapt the framework of [1] for the calculation of information bounds in semi-parametric models to the case where there are two unknown one-dimensional density functions. In Section 2.2 we specify the model $f$ for nested case control sampling and formally state our main result showing that the MPLE, and the Breslow estimator, achieve their respective efficiency lower bounds.

### 2.1. Information bounds in semi-parametric models

This section closely follows the treatment in [1] for deriving lower bounds for estimation in semi-parametric models; see also the text [3]. Let $L^2(\mu)$ denote the collection of functions that are square integrable with respect to a measure $\mu$, and for $u, v \in L^2(\mu)$ we let $\langle u, v \rangle_\mu = \int uv \, d\mu$ and $\|u\|_\mu^2 = \langle u, u \rangle_\mu$. Here, as in [1], the data consists of $n$ i.i.d. observations $X_1, \ldots, X_n$ taking values in a measurable space $(\mathcal{X}, \mathcal{F}_\mathcal{X})$, and the density function $f$ of a single observation is with respect to a sigma-finite measure $\sigma$. We consider a model where the density $f = f(\cdot, \theta, g, h)$ is determined by a real parameter $\theta$, the one of most interest, and by the infinite dimensional parameter $p = (g, h)$, a vector of two unknown densities $g$ and $h$, the baseline failure time density, and the marginal covariate density, respectively.

Let $\mathcal{D}^+$ and $\mathcal{D}$ denote the collection of densities with respect to Lebesgue measure $\nu^+$ and $\nu$ on $\mathbb{R}^+ = [0, \infty)$ and $\mathbb{R}$, respectively. We let the parameter space $\mathcal{G}$ for the unknown baseline failure density be

$$\mathcal{G} = \mathcal{D}^+.$$

To impose growth conditions on the covariates similar to the ones typically assumed, for a covariate density $\mathsf{h} \colon \mathbb{R} \to [0, \infty)$ and $\theta \in \mathbb{R}$ let

$$M_\mathsf{h}(\theta) = \int \mathsf{h} \, d\nu_\theta, \qquad \text{where } \frac{d\nu_\theta}{d\nu} = e^{\theta z}.$$

For some fixed $\xi > 0$ and $0 < \theta_\xi < \theta_\kappa$ we let the parameter space for the covariate density be

$$\mathcal{H} = \{\mathsf{h} \in \mathcal{D} : M_\mathsf{h}(\theta) < \infty \text{ for all } |\theta| < \theta_\kappa \text{ and } M_\mathsf{h}(\theta_\xi) + M_\mathsf{h}(-\theta_\xi) < \xi\}.$$

Hence, the parameter space $\mathcal{P}$ for the pair $p$ of unknown densities is given by

$$\mathcal{P} = \mathcal{G} \times \mathcal{H}.$$

Adopting slightly inconsistent notation for the sake of ease, we let $\theta_0$ denote the null parameter in $\mathbb{R}$, and henceforth, $g$ and $h$ the null parameters in $\mathcal{G}$ and $\mathcal{H}$, respectively; we label them also as $g_0$ and $h_0$ when convenient.

For $\tau \in \mathbb{R}$ let $\Theta(\tau)$ denote the collection of all real sequences $\{\theta_n\}_{n \geq 1}$ such that

$$|\sqrt{n}(\theta_n - \theta_0) - \tau| \to 0 \quad \text{as } n \to \infty \quad \text{and set} \quad \Theta = \bigcup \{\Theta(\tau) : \tau \in \mathbb{R}\}.$$

Let $\Pi_\theta = L^2(\nu^+) \times L^2(\nu_\theta)$ and for $\gamma = (\alpha, \beta) \in \Pi_\theta$ let $\|\gamma\|_{\Pi_\theta} = \max\{\|\alpha\|_{\nu^+}, \|\beta\|_{\nu_\theta}\}$, the product metric, and, with $p = (g, h)$ as the null parameter, let $\mathcal{C}(p, \gamma)$ be the collection of all sequences $\{p_n\}_{n \geq 0} = \{(g_n, h_n)\}_{n \geq 0} \subset \mathcal{P}$ such that

$$\|\sqrt{n}(p_n^{1/2} - p^{1/2}) - \gamma\|_{\Pi_\theta} \to 0 \quad \text{as } n \to \infty \text{ for all } |\theta| < \theta_\kappa. \tag{2}$$

Let $\Gamma$ be the set of all $\gamma$ such that (2) holds for some $\{p_n\}_{n \geq 0} \subset \mathcal{P}$, and

$$\mathcal{C}(p) = \bigcup_{\gamma \in \Gamma} C(p, \gamma).$$

By considering the components of $\{p_n\}_{n \geq 0}$ we see that $\mathcal{C}_1(g, \alpha)$ is the collection of all sequences $\{g_n\}_{n \geq 0}$ in $\mathcal{G}$, that satisfy

$$\|\sqrt{n}(g_n^{1/2} - g^{1/2}) - \alpha\|_{\nu^+} \to 0 \quad \text{as } n \to \infty,$$

and therefore $\alpha \in L^2(\nu^+)$ satisfies $\alpha \perp g^{1/2}$ in $L^2(\nu^+)$, that is, $\langle \alpha, g^{1/2} \rangle_{\nu^+} = 0$, or,

$$\int_0^\infty g^{1/2} \alpha \, d\nu^+ = 0. \tag{3}$$

Now let

$$\mathcal{A} = \{\alpha \in L^2(\nu^+) : \text{there exists} \{g_n\}_{n \geq 0} \subset \mathcal{G} \text{ such that } \|\sqrt{n}(g_n^{1/2} - g^{1/2}) - \alpha\|_{\nu^+} \to 0\}$$

and set

$$\mathcal{C}_1(g) = \bigcup_{\alpha \in \mathcal{A}} \mathcal{C}_1(g, \alpha).$$

Similarly, $\mathcal{C}_2(h, \beta)$ is the collection of all sequences $\{h_n\}_{n \geq 0}$ in $\mathcal{H}$ such that

$$\|\sqrt{n}(h_n^{1/2} - h^{1/2}) - \beta\|_{\nu_\theta} \to 0 \quad \text{as } n \to \infty \text{ for all } |\theta| < \theta_\kappa. \tag{4}$$

For $\theta = 0$ (4) yields

$$\|\sqrt{n}(h_n^{1/2} - h^{1/2}) - \beta\|_\nu \to 0 \qquad \text{as } n \to \infty, \tag{5}$$

and therefore that $\beta$ satisfies

$$\int_{-\infty}^{\infty} h^{1/2} \beta \, d\nu = 0. \tag{6}$$

Now let $\mathcal{B}$ be the collection of all $\beta \in L^2(\nu)$ such that there exists $\{h_n\}_{n \geq 0} \subset \mathcal{H}$ such that

$$\|\sqrt{n}(h_n^{1/2} - h^{1/2}) - \beta\|_{\nu_\theta} \to 0 \qquad \text{for all } |\theta| < \theta_\kappa,$$

and set

$$\mathcal{C}_2(h) = \bigcup_{\beta \in \mathcal{B}} \mathcal{C}_2(h, \beta).$$

Clearly

$$\mathcal{C}(p, \gamma) = \mathcal{C}_1(g, \alpha) \times \mathcal{C}_2(h, \beta), \qquad \mathcal{C}(p) = \mathcal{C}_1(g) \times \mathcal{C}_2(h) \quad \text{and} \quad \Gamma = \mathcal{A} \times \mathcal{B}.$$

The following three assumptions will be needed to demonstrate Theorems 2.1 and 2.2, and, in addition, the fourth will be needed for Theorems 2.3 and 2.4. The first is that $\Gamma$ is a subspace of $L^2(\nu^+) \times L^2(\nu)$ or, equivalently,

**Assumption 2.1.** *The sets $\mathcal{A}$ and $\mathcal{B}$ are subspaces of $L^2(\nu^+)$ and $L^2(\nu)$, respectively.*

It is shown in [1] that parts of the following assumption are a consequence of the Hellinger differentiability of $f$; we verify Assumption 2.2 directly.

**Assumption 2.2.** *There exists $\rho_\theta \in L^2(\sigma)$ and linear operators $A \colon L^2(\nu^+) \to L^2(\sigma)$ and $B \colon L^2(\nu) \to L^2(\sigma)$ such that for any $(\tau, \alpha, \beta) \in \mathbb{R} \times \mathcal{A} \times \mathcal{B}$ and*

$$(\{\theta_n\}_{n \geq 0}, \{g_n\}_{n \geq 0}, \{h_n\}_{n \geq 0}) \in \Theta(\tau) \times \mathcal{C}_1(g, \alpha) \times \mathcal{C}_2(h, \beta), \tag{7}$$

*the sequence of densities given by $f_n = f(\cdot, \theta_n, g_n, h_n)$ for $n = 0, 1, \ldots$ satisfies*

$$\|\sqrt{n}(f_n^{1/2} - f_0^{1/2}) - \zeta\|_\sigma \to 0 \qquad \text{for } \zeta = \tau \rho_\theta + A\alpha + B\beta \text{ as } n \to \infty. \tag{8}$$

Let

$$\mathbb{H} = \{\zeta \in L^2(\sigma) : \zeta = \tau \rho_\theta + A\alpha + B\beta \text{ for some } \tau \in \mathbb{R}, \alpha \in \mathcal{A}, \beta \in \mathcal{B}\} \tag{9}$$

and

$$\mathbb{K} = \{\delta \in L^2(\sigma) : \delta = A\alpha + B\beta \text{ for some } \alpha \in \mathcal{A} \text{ and } \beta \in \mathcal{B}\}. \tag{10}$$

The classical projection theorem shows that the orthogonal projection of $\rho_\theta$ onto the closure of $\mathbb{K}$ is an element of the closure of $\mathbb{K}$. However, we consider situations satisfying the following assumption, that is, where $\mathbb{K}$ itself contains the projection of $\rho_\theta$.

**Assumption 2.3.** *There exists $\hat{\alpha} \in \mathcal{A}$ and $\hat{\beta} \in \mathcal{B}$ such that $\hat{\delta} = A\hat{\alpha} + B\hat{\beta}$ satisfies*

$$\rho_\theta - \hat{\delta} \perp \delta \quad \text{for all } \delta \in \mathbb{K}.$$

Since for any $\delta = A\alpha + B\beta \in \mathbb{K}$, by orthogonality,

$$\begin{aligned}
\|\rho_\theta - \delta\|_\sigma^2 &= \|\rho_\theta - A\alpha - B\beta\|_\sigma^2 \\
&= \|\rho_\theta - \hat{\delta} - A(\alpha - \hat{\alpha}) - B(\beta - \hat{\beta})\|_\sigma^2 \\
&= \|\rho_\theta - \hat{\delta}\|_\sigma^2 + \|A(\alpha - \hat{\alpha}) + B(\beta - \hat{\beta})\|_\sigma^2 \\
&\geq \|\rho_\theta - \hat{\delta}\|_\sigma^2,
\end{aligned}$$

hence $\hat{\delta}$ minimizes $\|\rho_\theta - \delta\|_\sigma^2$ over $\delta \in \mathbb{K}$, and thus corresponds to the worst case direction of approach to the null, that is, the one that minimizes the available information. Set the effective information $I_*$ to be

$$I_* = 4\|\rho_\theta - \hat{\delta}\|_\sigma^2. \tag{11}$$

For $\zeta \in \mathbb{H}$, let $\mathcal{F}(f, \zeta)$ be the collection of all sequences $\{f_n\}_{n \geq 0}$ such that (8) holds, and $\mathcal{F}(f)$ the union of $\mathcal{F}(f, \zeta)$ over all $\zeta \in \mathbb{H}$. We say that an estimator $\hat{\theta}_n$ of $\theta_0$ is *regular* at $f = f(\cdot, \theta_0, g, h)$ if for every sequence $f_n(\cdot, \theta_n, g_n, h_n)$ with $\{\theta_n\}_{n \geq 0}, \{g_n\}_{n \geq 0}$ and $\{h_n\}_{n \geq 0}$ as in (7), the distribution of $\sqrt{n}(\hat{\theta}_n - \theta_0)$ converges in distribution to $\mathcal{L} = \mathcal{L}(f)$, which depends on $f$ but not on the particular sequence $f_n$.

The setup above differs in two ways from that in [1]. First, the model considered here has two nonparametric components, $g$ and $h$, while in [1] only one nonparametric component is considered. Second, as we specify the parameter space $\mathcal{H}$ on the covariate density $h$ in such a way as to accommodate more relaxed integrability conditions, the resulting space of perturbations $\mathcal{B}$ is expressed as the intersection of subspaces (see [9]), one for each $\theta$ in $(-\theta_\kappa, \theta_\kappa)$. This is so as the perturbations $\beta$ are required to be limiting approximations to $\sqrt{n}(h_n^{1/2} - h^{1/2})$ in $L^2(\nu_\theta)$ for all $|\theta| < \theta_\kappa$, rather than in $L^2(\nu)$. As (4) implies (5) our condition gives rise to a smaller collection $\mathcal{B}$ of perturbations than in [1]. Nevertheless, only minimal adaptations of the proofs of Theorems 3.1 and 3.2 and Theorems 4.1 and 4.2 of [1] are required to demonstrate Theorems 2.1–2.4 for our model, so these are relegated to Section A.4.

**Theorem 2.1.** *Suppose that $\hat{\theta}_n$ is a regular estimator of $\theta_0$ in the model $f = f(\cdot, \theta, g, h)$ with limit law $\mathcal{L} = \mathcal{L}(f)$ and that assumptions 2.1–2.3 hold. Then $\mathcal{L}$ is the convolution of a normal $\mathcal{N}(0, 1/I_*)$ distribution with a distribution depending only on $f$, where $I_*$ is given by (11).*

We may also adapt the asymptotic minimax result of [1]. Recall that we say a loss function $\ell:\mathbb{R} \to \mathbb{R}^+$ is subconvex when $\{x:\ell(x) \leq y\}$ is closed, convex and symmetric for every $y \geq 0$. We will also assume our loss function satisfies

$$\int_{-\infty}^{\infty} \ell(z)\phi(az)\,\mathrm{d}z < \infty \qquad \text{for all } a > 0, \tag{12}$$

where $\phi$ denotes the standard normal density function.

**Theorem 2.2.** *Suppose Assumptions 2.1–2.3 hold and that $\ell$ is subconvex and satisfies (12). For $c \geq 0$ let*

$$B_n(c) = \{f_n \in \mathcal{F} : \sqrt{n}\|f_n^{1/2} - f^{1/2}\|_\sigma \leq c\}. \tag{13}$$

*Then*

$$\lim_{c \to \infty} \liminf_{n \to \infty} \inf_{\hat{\theta}_n} \sup_{f_n \in B_n(c)} E_{f_n}\ell(\sqrt{n}(\hat{\theta}_n - \theta_n)) \geq E\ell(Z_*), \tag{14}$$

*where $Z_* \sim \mathcal{N}(0, 1/I_*)$ and $I_*$ is given by (11).*

The infimum in (14) is taken over the class of "generalized procedures," the closure of the class of randomized Markov kernel procedures (see [13], page 235). We also obtain lower bounds on the performance of regular estimators of the baseline survival function $\overline{G}(\cdot)$ by similarly adapting Theorems 4.1 and 4.2 of [1] under the following assumption.

**Assumption 2.4.** *The linear operator $A^*A : L^2(\nu^+) \to L^2(\nu^+)$ is invertible with bounded inverse $(A^*A)^{-1}$.*

We also suppose that, perhaps by a suitable map such as the probability integral transformation, the density $g$ is supported on $[0,1]$. Let

$$G_s = (I_{[0,s]} - G(s))g(s)^{1/2},$$

and define the covariance functions

$$K(s,t) = \langle G_s, (A^*A)^{-1} G_t \rangle_{\nu^+} \quad \text{and} \quad K_*(s,t) = K(s,t) + 4I_*^{-1} \int_0^s \hat{\alpha} g^{1/2} \int_0^t \hat{\alpha} g^{1/2}, \tag{15}$$

where $I_*$ is given by (11) and $\hat{\alpha}$ is as in Assumption 2.3. For the precise definition of a regular estimator of $G(\cdot)$, analogous to that for estimators of $\theta_0$, see [1].

**Theorem 2.3.** *Suppose that $\widehat{G(\cdot)}_n$ is a regular estimator of $G(\cdot) = \int_0^\cdot g\,\mathrm{d}\nu^+$ in the model $f = f(\cdot;\theta,g,h)$ with limit process $\mathbb{S}$, that Assumptions 2.2–2.4 hold, and that Assumption 2.1 holds with $\mathcal{A}$ given by $\{\alpha \in L^2(\nu^+) : \int \alpha g^{1/2}\,\mathrm{d}\nu^+ = 0\}$. Then*

$$\mathbb{S} =_d \mathbb{Z}_* + \mathbb{W},$$

where $\mathbb{Z}_*$ *is a mean zero Gaussian process with covariance function* $K_*(s,t)$ *given by (15) and the process* $\mathbb{W}$ *is independent of* $\mathbb{Z}_*$.

For the local asymptotic minimax bound, we let $\ell\colon C[0,1]\to\mathbb{R}^+$ be a subconvex loss function, such as $\ell(x)=\sup_t|x(t)|, \ell(x)=\int|x(t)|^2\,dt$, or $\ell(x)=\mathbf{1}(x\colon\|x\|\geq c)$.

**Theorem 2.4.** *Suppose the hypotheses of Theorem 2.3 are satisfied, that $\ell$ is subconvex, and that $B_n(c)$ is as in (13). Then*

$$\lim_{c\to\infty}\lim_{n\to\infty}\inf_{\widehat{G(\cdot)}_n}\sup_{f_n\in B_n(c)} E_{f_n}\ell(\sqrt{n}(\widehat{G(\cdot)}_n - G_n)) \geq E\ell(\mathbb{Z}_*), \tag{16}$$

*where $\mathbb{Z}_*$ is the mean zero Gaussian process with covariance $K_*(s,t)$ given by (15).*

The infimum over estimators $\widehat{G(\cdot)}_n$ in (16) is taken over the class of "generalized procedures" as in [13], page 235. The proofs of Theorems 2.1–2.4 in the Appendix detail the modifications required for the application of the methods of [1] to the case at hand.

## 2.2. Main results

We now specify our model $f$ for the nested case control sampling of $m-1$ controls for the failure in each group. For any integer $k$, let $[k]=\{1,\ldots,k\}$, and for any set $S$ let $\mathcal{P}_k(S)$ be the collection of all subsets of $S$ of size $k$. Groups of individuals of size $\eta\geq m$ are observed up to the time of the first failure, at which point covariates are collected on a simple random sample of $m-1$ non-failed individuals and the failure.

An observation $X=(\eta,i,r,t,z_r)$ consists of the group size $\eta$, the identity $i\in[\eta]$ of the failed individual, the group $r\subset[\eta]$ of the $m$ individuals whose covariates are collected, the time $t$ of the failure, and the covariates $z_r=\{z_j,j\in r\}$. In particular, $X$ takes values in the space

$$\mathcal{X}=\bigcup_{\eta\geq m}\{\eta\times[\eta]\times\mathcal{P}_m([\eta])\times\mathbb{R}^+\times\mathbb{R}^m\},$$

which we endow with the $\sigma$-finite product measure

$$\sigma=(\text{counting measure})\times(\text{counting measure})\times(\text{counting measure})\times\nu^+\times\nu^m.$$

To begin the specification of the density $f$ of the observations, corresponding to the baseline survival density $g$ on $\mathbb{R}^+$ are the baseline survival and hazard functions, for $t\geq 0$, given by, respectively

$$\overline{G}(t)=\int_t^\infty g(u)\,du \quad\text{and}\quad \lambda(t)=\begin{cases} g(t)\overline{G}^{-1}(t), & \text{for } \overline{G}(t)>0, \\ 0, & \text{otherwise.} \end{cases}$$

Under the assumed standard exponential relative risk form, the hazard function $\lambda(t;z)$ for an individual with covariate value $z$ is the baseline hazard scaled by the factor $\exp(\theta z)$, that is, $\lambda(t;z) = \exp(\theta z)\lambda(t)$, resulting in survival and density functions, respectively, of

$$\overline{G}_\theta(t;z) = \overline{G}^{e^{\theta z}}(t) \quad \text{and} \quad g_\theta(t;z) = \begin{cases} e^{\theta z} g(t) \overline{G}^{e^{\theta z}-1}(t), & \text{for } \overline{G}(t) > 0, \\ 0, & \text{otherwise;} \end{cases}$$

we note $g_0(t;z) = g_\theta(t;0) = g(t)$. As the marginal covariate density is $h$, the survival function $\overline{G}_\theta(t;z)$ averaged over individuals with covariate density $h(z)$ results in the (mixture) survival function

$$\overline{G}_\theta(t) = \int \overline{G}_\theta(t;z) h(z) \, dz$$

for individuals whose covariates are not observed.

The group size $\eta$ may vary from strata to strata, and we assume it to be random with distribution, say, $\varrho$. At the time $t$ of the failure of individual $i$, a simple random sample of size $m-1$ is taken from the non-failures to serve as controls. Hence, when the group size is $\eta$ and the identity of the failure $i$, the probability that the set $r \subset [\eta]$ is selected is given by

$$K_{\eta,m} = \binom{\eta-1}{m-1}^{-1}$$

for any set $r$ of size $m$ containing $i$. We assume that the individuals in $[\eta]$ are independent, and therefore the density of the sampled covariates $z_r$ is the product

$$h(z_r) = \prod_{j \in r} h(z_j).$$

Putting all the factors together, the density for $X = (\eta, i, r, t, z_r)$ is given by

$$f(X;\theta,g,h) = K_{\eta,m} e^{\theta z_i} g(t) \overline{G}(t)^{\sum_{j \in r} e^{\theta z_j} - 1} \overline{G}_\theta(t)^{\eta-m} h(z_r) \varrho(\eta) \quad (17)$$

$$= K_{\eta,m} g(t;z_i) \left[ \prod_{j \in r \setminus \{i\}} \overline{G}(t;z_j) \right] \overline{G}_\theta(t)^{\eta-m} h(z_r) \varrho(\eta).$$

For the sake of clarity or brevity, the density may be written with either its parameters or its variables suppressed, that is, as $f(\eta, i, r, t, z_r)$ or $f(\theta, g, h)$, respectively. At the null, (17) reduces to

$$f(X;\theta_0, g, h) = K_{\eta,m} g(t) \overline{G}(t)^{\eta-1} h(z_r) \varrho(\eta), \quad (18)$$

which, in agreement with the notation introduced in Section 2.1, may appear in the abbreviated form $f_0$. We may take the distribution $\varrho$ of $\eta$ as known when proving Theorem 2.5 since the MPLE is computed without knowledge of $\varrho$ and already achieves the bound (20) in the limit.

We are now ready to state our main result regarding the estimation of the parametric component of the model.

**Theorem 2.5.** *Suppose that $\eta \geq 2$ almost surely, $E[\eta^5] < \infty$ and at least one of the following conditions is satisfied:*

(i) *Positivity: The parameter space $\Theta = [0, \infty)$, the covariates $Z$ take on non-negative values, and $\eta \geq m$ almost surely.*

(ii) *Boundedness: The covariates $Z$ are bounded and $\eta \geq m$ almost surely.*

(iii) *Cohort size: $1 \leq m \leq \eta - 4$ almost surely.*

*Then Theorems 2.1 and 2.2 obtain for the nested case control model given in (17) with effective information*

$$I_*^\varrho(\theta_0) = \mathrm{Var}(Z)\left(1 - \frac{1}{m}\right) + m\,\mathrm{Var}(Z)\left(2\,\mathrm{Var}\left(\frac{1}{\eta}\right) + \left[E\left(\frac{1}{\eta}\right)\right]^2\right). \quad (19)$$

*In particular, under any of the above three scenarios, if $\varrho_n$ is a sequence of distributions such that $\eta_n \to_p \infty$ when $\eta_n$ has distribution $\varrho_n$, then*

$$I_*(\theta_0) = \lim_{n \to \infty} I_*^{\varrho_n}(\theta_0) = \mathrm{Var}(Z)\left(\frac{m-1}{m}\right) \quad (20)$$

*and hence the Cox MPLE is efficient for the limiting nested case control model.*

The situation where there is full cohort information is covered by the special case $P(\eta = m) = 1$, for which (19) reduces to the lower bound $\mathrm{Var}(Z)$, recovering the result of [1] for the case of no censoring. See Section A.3 for some remarks on the rationale behind the three conditions in Theorem 2.5.

Next, we consider lower bounds for the estimation of the nonparametric component of the model. It is shown in [8] that the Breslow estimator of the baseline survival is asymptotically normal with covariance function

$$\omega(s,t) = \overline{G}(t)\overline{G}(s)\left(\int_0^{s \wedge t} \frac{\mathrm{d}G}{E[\eta\overline{G}(u)^{\eta+1}]} + [E(Z)]^2(\log\overline{G}(t)\log\overline{G}(s))[I_*(\theta_0)]^{-1}\right), \quad (21)$$

where $I_*(\theta_0)$ is given in (20).

**Theorem 2.6.** *Let the hypotheses of Theorem 2.5 be satisfied. Then on any interval $[0, T_0]$ for which $\overline{G}(T_0) > 0$, the conclusions of Theorems 2.3 and 2.4 hold with*

$$K_*(s,t) = \overline{G}(t)\overline{G}(s)\left(\int_0^{s \wedge t} \frac{\mathrm{d}G}{E[\eta\overline{G}(u)^{\eta+1}]} + [E(Z)]^2(\log\overline{G}(t)\log\overline{G}(s))[I_*^\varrho(\theta_0)]^{-1}\right). \quad (22)$$

By (20) and (21), we see that the Breslow estimator becomes asymptotically efficient as the cohort size increases under the nested case control model considered.

Theorem 2.5 follows from Theorems 2.1 and 2.2. The application of these theorems is a consequence of Theorem 4.1, which provides the effective information $I_*^\varrho(\theta_0)$, and the verification of Assumptions 2.1–2.3. In [9], a simple argument shows that Assumption 2.1 is satisfied with $\mathcal{A}$ and $\mathcal{B}$ given by (40). The verification of Assumption 2.2 is somewhat involved. The relevant quantities, $A, B, \hat{\alpha}, \hat{\beta}$ and $\rho_0$, are given in (24), (25), Lemmas 3.2, 3.4 and (23), respectively. The remainder of the verification of Assumption 2.2, that is, the convergence to zero in (8), is shown in Lemma 3.1 whose proof is deferred to [9]. Assumption 2.3 follows in a fairly straightforward manner from (40). Some remarks on the calculations in [9] can be found in Section A.3.

Theorem 2.6 follows similarly from Theorems 2.3 and 2.4. In addition to Assumptions 2.1–2.3, the application of these theorems follow from Theorem 4.2, which verifies the covariance lower bound (22); Lemma A.2, from which Assumption 2.4 on $[0, T_0]$ follows easily; and (40), which shows that $\mathcal{A}$ is of the form required by Theorem 2.3. Regarding the restriction of the result to $[0, T_0]$, see example 4 in [1], page 450 in particular, and the proof of Lemma 2 in [16].

## 3. Operators $A$ and $B$: properties

The following lemma provides the parametric score $\rho_0$ and the operators $A$ and $B$ required by Assumption 2.2 and needed for the computation of the effective information $I_*$ in (11). Sums over $r$ denote a sum over all $r \subset [\eta]$ of size $m$, and sums over $\eta, i, r$ are short for the sum over all $\eta \in \mathbb{Z}^+$, $i \in [\eta]$ and $r \subset [\eta]$ of size $m$ with $r \ni i$.

**Lemma 3.1.** *Assumption 2.2 is satisfied for the nested case control model (17) with*

$$\rho_0 = \frac{1}{2}\left[z_i + \log \overline{G}(t) \sum_{j \in r}(z_j - EZ) + \eta EZ \log \overline{G}(t)\right] f_0^{1/2}, \tag{23}$$

$$A\alpha = \left(g^{-1/2}(t)\alpha(t) + \frac{(\eta-1)\int_t^\infty g^{1/2}\alpha\,d\nu}{\overline{G}(t)}\right) f_0^{1/2} \tag{24}$$

*and*

$$B\beta = \left(\sum_{j \in r} h^{-1/2}(z_j)\beta(z_j)\right) f_0^{1/2}. \tag{25}$$

Lemma 3.1 is proved in [9].

### 3.1. $A$ operator: properties

Regarding the definition and calculation of adjoint operators such as $A^*$ in the following lemma, the reader is referred to [12]. The proof of the following lemma appears in Section A.1.

**Lemma 3.2.** *Let $\rho_0$ and $A$ be given by (23) and (24), respectively. Then the function*

$$\hat{\alpha} = \frac{EZ}{2}[1 + \log \overline{G}(t)]g^{1/2}(t)$$

*is the solution to the normal equation $A^*A\alpha = A^*\rho_0$ and the projection of $\rho_0$ onto the range of $A$ is given by*

$$A\hat{\alpha} = \frac{EZ}{2}[1 + \eta \log \overline{G}(t)]f_0^{1/2}. \tag{26}$$

### 3.2. $B$ operator: properties

Let $r \subset [\eta]$ of size $m$ be fixed. For $s \subset r$ let $z_s = \{z_j : j \in s\}$ and $z_{\neg s} = \{z_j : j \in r \setminus s\}$ and denote integration over $z_s$ and $z_{\neg s}$ with respect to the measures $\nu^{|s|}$ and $\nu^{m-|s|}$ by $\mathrm{d}z_s$ and $\mathrm{d}z_{\neg s}$, respectively. When $s = \{j\}$, we identify that $j$th variable $z_j$ with $z$. Integration with respect to $\nu^+$ is often indicated by $\mathrm{d}t$, but may also be indicated by other notations such as $\mathrm{d}u$, or suppressed, when clear from context.

**Lemma 3.3.** *The adjoint $B^* : L^2(\sigma) \to L^2(\nu)$ of the operator $B$ in (25) is given by*

$$B^*\mu = h^{-1/2}(z) \sum_{\eta,i,r,j \in r} \int_{z_{\neg j}} \int_0^\infty f_0^{1/2} \mu \, \mathrm{d}t \, \mathrm{d}z_{\neg j}. \tag{27}$$

**Proof.** As $B = \sum_{j \in r} B_j$ with

$$B_j\beta = h^{-1/2}(z_j)f_0^{1/2}\beta(z_j) \qquad \text{for } \beta \in L^2(\nu),$$

by linearity one need only sum the adjoints $B_j^*$ of $B_j$ over $j \in r$ to obtain $B^*$. For $\mu \in L^2(\sigma)$, the calculation

$$\begin{aligned}
\langle B_j\beta, \mu \rangle_\sigma &= \int_{\mathcal{X}} B_j\beta\mu \, \mathrm{d}\sigma \\
&= \sum_{\eta,i,r} \int_{z_r} \int_0^\infty h^{-1/2}(z_j)f_0^{1/2}\beta(z_j)\mu \, \mathrm{d}t \, \mathrm{d}z_r \\
&= \int_z \beta(z)\left(h^{-1/2}(z) \sum_{\eta,i,r} \int_{z_{\neg j}} \int_0^\infty f_0^{1/2}\mu \, \mathrm{d}t \, \mathrm{d}z_{\neg j}\right) \mathrm{d}z \\
&= \langle \beta, B_j^*\mu \rangle_\nu
\end{aligned}$$

provides the desired conclusion. $\square$

The proof of the following lemma appears in Section A.2.

**Lemma 3.4.** *The function*

$$\hat{\beta} = \frac{1}{2} h^{1/2}(z) E\left[\frac{\eta - m}{m\eta}\right](z - EZ) \tag{28}$$

*is the solution to the normal equation $B^*B\beta = B^*\rho_0$, and the projection of $\rho_0$ onto the range of $B$ is given by*

$$B\hat{\beta} = \frac{1}{2}\left(E\left[\frac{\eta - m}{m\eta}\right]\sum_{j \in r}(z_j - EZ)\right) f_0^{1/2}. \tag{29}$$

## 4. Lower bound calculations

We begin the computation of the information bound by showing that the two operators $A$ and $B$ have orthogonal ranges.

**Lemma 4.1.** *Let $A$ and $B$ be the operators given by (24) and (25), respectively. Then*

$$B^*A = 0 \quad \text{and} \quad A^*B = 0.$$

**Proof.** Since $(A^*B)^* = B^*A$ it suffices to prove only the first claim. By (24) and (27),

$$B^*A\alpha = B^*\left(g^{-1/2}(t)\alpha(t) + \frac{(\eta - 1)\int_t^\infty g^{1/2}\alpha}{\overline{G}(t)}\right) f_0^{1/2}$$

$$= h^{-1/2}(z) \sum_{\eta, i, r, j \in r} \int_{z_{\neg j}} \int_0^\infty \left(g^{-1/2}(t)\alpha(t) + \frac{(\eta - 1)\int_t^\infty g^{1/2}\alpha}{\overline{G}(t)}\right) f_0 \, dt \, dz_{\neg j}$$

$$= h^{-1/2}(z) \sum_\eta K_{\eta,m}\varrho(\eta)$$

$$\times \left[\sum_{i, r, j \in r} \int_{z_{\neg j}} h(z_r) \int_0^\infty \left(g^{-1/2}(t)\alpha(t)\right.\right.$$

$$\left.\left. + \frac{(\eta - 1)\int_t^\infty g^{1/2}\alpha}{\overline{G}(t)}\right) g(t)\overline{G}^{\eta-1}(t) \, dt \, dz_{\neg j}\right].$$

Integrating the inner integral by parts,

$$\int_0^\infty \left(g^{-1/2}(t)\alpha(t) + (\eta - 1)\frac{\int_t^\infty g^{1/2}\alpha}{\overline{G}(t)}\right) g(t)\overline{G}^{\eta-1}(t) \, dt$$

$$= \int_0^\infty \left( g^{1/2}(t) \overline{G}^{\eta-1}(t) \alpha(t) + (\eta-1) g(t) \overline{G}^{\eta-2}(t) \int_t^\infty g^{1/2} \alpha \right) \mathrm{d}t$$

$$= \int_0^\infty g^{1/2}(t) \overline{G}^{\eta-1}(t) \alpha(t) \, \mathrm{d}t - \int_0^\infty (\overline{G}^{\eta-1}(t))' \left( \int_t^\infty g^{1/2} \alpha \right) \mathrm{d}t$$

$$= \int_0^\infty g^{1/2}(t) \overline{G}^{\eta-1}(t) \alpha(t) \, \mathrm{d}t - \left[ \overline{G}^{\eta-1}(t) \int_t^\infty g^{1/2} \alpha \right]_0^\infty - \int_0^\infty g^{1/2}(t) \overline{G}^{\eta-1}(t) \alpha(t) \, \mathrm{d}t$$

$$= \int_0^\infty g^{1/2} \alpha,$$

which equals zero by (3). □

The perpendicularity relation that holds between $A$ and $B$ allows for the application of the following lemma, which simplifies the calculation of the information bound.

**Lemma 4.2.** *Let $\mathbb{K}$ be given by (9). Then under the perpendicularity relations provided by Lemma 4.1, the function*

$$\hat{\delta} = A\hat{\alpha} + B\hat{\beta} \qquad \textit{minimizes } \|\rho_0 - \delta\|_\sigma \textit{ over } \delta \in \mathbb{K},$$

*where $\hat{\alpha}$ and $\hat{\beta}$ are the solutions to the normal equations $A^*A\alpha = A^*\rho_0$ and $B^*B\beta = B^*\rho_0$, respectively. Consequently, the effective information (11) is given by*

$$I_*(\theta_0) = 4\|\rho_0 - A\hat{\alpha} - B\hat{\beta}\|_\sigma^2. \tag{30}$$

**Proof.** Since $A^*B = 0$ we have $A^*\rho_0 = A^*A\hat{\alpha} = A^*\hat{\delta}$ and similarly $B^*\rho_0 = B^*B\hat{\beta} = B^*\hat{\delta}$. Therefore $(A+B)^*\rho_0 = (A+B)^*\hat{\delta}$, or $(A+B)^*(\rho_0 - \hat{\delta}) = 0$. Hence we have

$$\rho_0 - \hat{\delta} \perp \mathbb{K} \quad \text{and} \quad \hat{\delta} \in \mathbb{K},$$

showing $\hat{\delta}$ is the claimed minimizer. □

We pause to record a simple calculation that will be used frequently in what follows.

**Lemma 4.3.** *Let $s(t)$ be any density on $\mathbb{R}^+$ and $\overline{S}(t)$ the corresponding survival function. Then for all integers $\eta$ and $k$ satisfying $\eta \geq k$, and $j = 1, 2, \ldots$,*

$$\int_0^\infty s(t) \overline{S}(t)^{\eta-k} [\log \overline{S}(t)]^j \, \mathrm{d}t = (-1)^j (\eta - k + 1)^{-(j+1)} j!.$$

*In particular, as $\log S(t) \leq 0$ for all $t \in \mathbb{R}^+$, if $k$ and $j$ are fixed, then for any constant $C > 1$ there exists $\eta_C$ such that*

$$\int_0^\infty s(t) \overline{S}(t)^{\eta-k} |\log \overline{S}(t)|^j \, \mathrm{d}t \leq \frac{Cj!}{\eta^{j+1}} \qquad \textit{for all } \eta \geq \eta_C.$$

**Proof.** Rewriting the integral and then applying the change of variables $u = \overline{S}(t)^{\eta-k+1}$ followed by $u = e^{-x}$ we have

$$\int_0^\infty s(t)\overline{S}(t)^{\eta-k}[\log \overline{S}(t)]^j \, dt = (\eta-k+1)^{-j} \int_0^\infty s(t)\overline{S}(t)^{\eta-k}[\log \overline{S}(t)^{\eta-k+1}]^j \, dt$$

$$= (\eta-k+1)^{-(j+1)} \int_0^1 [\log u]^j \, du$$

$$= (-1)^j(\eta-k+1)^{-(j+1)}\Gamma(j+1)$$

$$= (-1)^j(\eta-k+1)^{-(j+1)} j!.$$

Taking absolute value and noting that $(\eta-k+1)/\eta \to 1$ suffices to prove the final claim. $\square$

**Theorem 4.1.** *The effective information for the nested case control model (17) is given by (19).*

**Proof.** Substituting (23), (26) and (29) into (30) we obtain

$$I_*^\varrho(\theta_0) = \left\| \left[ (z_i - EZ) + \log \overline{G}(t) \sum_{j \in r}(z_j - EZ) - E\left[\frac{\eta-m}{\eta m}\right] \sum_{j \in r}(z_j - EZ) \right] f_0^{1/2} \right\|_\sigma^2$$

$$= \left\| \left[ (z_i - EZ) + \sum_{j \in r}(z_j - EZ)\left(\log \overline{G}(t) - E\left[\frac{\eta-m}{\eta m}\right]\right) \right] f_0^{1/2} \right\|_\sigma^2$$

$$= \left\| \left[ \left(1 + \log \overline{G}(t) - E\left[\frac{\eta-m}{\eta m}\right]\right)(z_i - EZ) \right.\right.$$

$$\left.\left. + \sum_{j \in r \setminus \{i\}}(z_j - EZ)\left(\log \overline{G}(t) - E\left[\frac{\eta-m}{\eta m}\right]\right) \right] f_0^{1/2} \right\|_\sigma^2,$$

which, by the independence of $Z_i$ and $\{Z_j, j \in r \setminus \{i\}\}$, equals

$$\left\| \left[ \left(1 + \log \overline{G}(t) - E\left[\frac{\eta-m}{\eta m}\right]\right)(z_i - EZ) \right] f_0^{1/2} \right\|_\sigma^2$$

$$+ \left\| \left[ \sum_{j \in r \setminus \{i\}}(z_j - EZ)\left(\log \overline{G}(t) - E\left[\frac{\eta-m}{\eta m}\right]\right) \right] f_0^{1/2} \right\|_\sigma^2.$$

Squaring and integrating against the null density (18) we obtain

$$\int_{\mathcal{X}} \left[ \left(1 + \log \overline{G}(t) - E\left[\frac{\eta-m}{\eta m}\right]\right)^2 (z_i - EZ)^2 \right.$$

$$+ \sum_{j \in r \setminus \{i\}} \left( \log \overline{G}(t) - E\left[\frac{\eta - m}{\eta m}\right] \right)^2 (z_j - EZ)^2 \Bigg] f_0 \, d\sigma$$

$$= \mathrm{Var}(Z) \sum_\eta K_{\eta,m} \varrho(\eta) \Bigg[ \sum_{i,r} \int_0^\infty \Bigg[ \left(1 + \log \overline{G}(t) - E\left[\frac{\eta - m}{\eta m}\right] \right)^2$$

$$+ (m-1) \left( \log \overline{G}(t) - E\left[\frac{\eta - m}{\eta m}\right] \right)^2 \Bigg] g(t) \overline{G}^{\eta-1}(t) \, dt \Bigg]$$

$$= \mathrm{Var}(Z) \sum_\eta \varrho(\eta) \Bigg[ \int_0^\infty \Bigg[ \left(1 + \log \overline{G}(t) - E\left[\frac{\eta - m}{\eta m}\right] \right)^2$$

$$+ (m-1) \left( \log \overline{G}(t) - E\left[\frac{\eta - m}{\eta m}\right] \right)^2 \Bigg] \eta g(t) \overline{G}^{\eta-1}(t) \, dt \Bigg]$$

$$= \mathrm{Var}(Z) \sum_\eta \varrho(\eta) \Bigg[ \int_0^\infty \Bigg[ 1 + 2 \left( \log \overline{G}(t) - E\left[\frac{\eta - m}{\eta m}\right] \right)$$

$$+ m \left( \log \overline{G}(t) - E\left[\frac{\eta - m}{\eta m}\right] \right)^2 \Bigg] \eta g(t) \overline{G}^{\eta-1}(t) \, dt \Bigg]$$

$$= \mathrm{Var}(Z) E \Bigg[ \left( 1 - 2 \left( \frac{1}{\eta} + E\left[\frac{\eta - m}{\eta m}\right] \right) \right.$$

$$+ m \left( \frac{2}{\eta^2} + 2 \frac{1}{\eta} E\left[\frac{\eta - m}{\eta m}\right] + \left( E\left[\frac{\eta - m}{\eta m}\right] \right)^2 \right) \Bigg],$$

by applying Lemma 4.3. Simplifying we obtain

$$I_*^\varrho(\theta_0) = \mathrm{Var}(Z) \left(1 - \frac{1}{m}\right) + m \, \mathrm{Var}(Z) \left( 2 \, \mathrm{Var}\left(\frac{1}{\eta}\right) + \left[E\left(\frac{1}{\eta}\right)\right]^2 \right),$$

which is (19). □

We now calculate the lower bound for the estimation of the baseline survival.

**Theorem 4.2.** *The covariance function $K_*(s,t)$ in (15) specializes to (22) for the nested case control model (17).*

**Proof.** Lemma A.2 shows that $A^*A$ is given by (36) with $M_0(t) = E[\eta \overline{G}^\eta(t)]$. Now (6.8) of [1] yields

$$K(s,t) = \overline{G}(t) \overline{G}(s) \int_0^{s \wedge t} \frac{dG}{M_0(u) \overline{G}(u)} = \overline{G}(t) \overline{G}(s) \int_0^{s \wedge t} \frac{dG}{E[\eta \overline{G}(u)^{\eta+1}]}.$$

Regarding the integral in (15), using the form $\hat{\alpha}$ given in Lemma 3.2, we have

$$\int_0^t \hat{\alpha} g^{1/2} \, d\nu = \int_0^t \frac{EZ}{2}[1 + \log \overline{G}(u)]g(u) \, du$$

$$= \frac{EZ}{2} \int_{\overline{G}(t)}^1 (1 + \log x) \, dx = -\frac{EZ}{2}\overline{G}(t) \log \overline{G}(t).$$

Substitution into (15) now yields (22). □

# Appendix

In the following four sections of this Appendix we prove Lemmas 3.2 and 3.4 and provide some remarks regarding the verification of Assumptions 2.1–2.3, and the proofs of Theorems 2.1–2.4.

## A.1. *A* operator

In this section we provide the proof to Lemma 3.2. We begin by calculating the adjoint $A^* \colon L^2(\sigma) \to L^2(\nu^+)$ of the operator $A$ given in Lemma 3.1.

**Lemma A.1.** *For the operator $A \colon L^2(\nu^+) \to L^2(\sigma)$ given in (24), write*

$$A = A_1 + A_2,$$

*where*

$$A_1 \alpha = g^{-1/2} f_0^{1/2} \alpha \quad \text{and} \quad A_2 \alpha = \frac{(\eta - 1) \int_t^\infty g^{1/2} \alpha}{\overline{G}(t)} f_0^{1/2}. \tag{31}$$

*Then the adjoint of $A$ is given by $A^* = A_1^* + A_2^*$, where*

$$A_1^* \mu = g^{-1/2}(t) \sum_{\eta, i, r} \int_{z_r} f_0^{1/2} \mu \, dz_r \quad \text{and}$$

$$A_2^* \mu = g^{1/2}(t) \sum_{\eta, i, r} (\eta - 1) \int_0^t \int_{z_r} \frac{f_0^{1/2}}{\overline{G}(u)} \mu \, dz_r \, du. \tag{32}$$

**Proof.** Let $\alpha \in L^2(\nu^+)$ and $\mu \in L^2(\sigma)$. Then

$$\langle A_1 \alpha, \mu \rangle_\sigma = \langle g^{-1/2} f_0^{1/2} \alpha, \mu \rangle_\sigma$$

$$= \sum_{\eta, i, r} \int_0^\infty \int_{z_r} g^{-1/2}(t) \alpha(t) f_0^{1/2} \mu \, dz_r \, dt$$

$$= \int_0^\infty \alpha(t)\left(g^{-1/2}(t)\sum_{\eta,i,r}\int_{z_r} f_0^{1/2}\mu\,\mathrm{d}z_r\right)\mathrm{d}t$$

$$= \langle \alpha, A_1^*\mu\rangle_{\nu^+}$$

when $A_1^*$ is as given in (32).

Next, writing $A_2$ as

$$A_2\alpha = L\int_t^\infty g^{1/2}\alpha \qquad \text{for } L = (\eta-1)\overline{G}^{-1}(t)f_0^{1/2}, \tag{33}$$

we have

$$\langle A_2\alpha,\mu\rangle_\sigma = \left\langle L\int_t^\infty g^{1/2}\alpha,\mu\right\rangle_\sigma$$

$$= \sum_{\eta,i,r}\int_0^\infty\int_{z_r}\int_t^\infty Lg^{1/2}\alpha\mu\,\mathrm{d}u\,\mathrm{d}z_r\,\mathrm{d}t$$

$$= \int_0^\infty \alpha(t)\left(g^{1/2}(t)\sum_{\eta,i,r}\int_0^t\int_{z_r} L\mu\,\mathrm{d}z_r\,\mathrm{d}u\right)\mathrm{d}t$$

$$= \langle\alpha, A_2^*\mu\rangle_{\nu^+}$$

when

$$A_2^*\mu = g^{1/2}(t)\sum_{\eta,i,r}\int_0^t\int_{z_r} L\mu\,\mathrm{d}z_r\,\mathrm{d}u.$$

Substituting $L$ from (33) now yields the stated conclusion. □

To help express the solution to the normal equations in $A$, for $\alpha \in L^2(\nu^+)$ define the operator $R$ as in [1] by the first equality in

$$R\alpha = g^{-1/2}(t)\alpha(t) - \frac{\int_t^\infty g^{1/2}\alpha}{\overline{G}(t)} = g^{-1/2}(t)\alpha(t) + \frac{\int_0^t g^{1/2}\alpha}{\overline{G}(t)}; \tag{34}$$

the second equality follows from (3). Also, set

$$M_0(t) = E[\eta\overline{G}^\eta(t)] \quad \text{and} \quad M_1(t) = E[Z]E[\eta\overline{G}(t)^\eta]. \tag{35}$$

**Lemma A.2.** *Let the operator $A$ be given by (24). Then, for $\alpha \in L^2(\nu^+)$,*

$$A^*A\alpha = \left[R\alpha(t)\frac{M_0(t)}{\overline{G}(t)} - \int_0^t R\alpha\frac{M_0}{\overline{G}(u)}\frac{\mathrm{d}G}{\overline{G}(u)}\right]g^{1/2}, \tag{36}$$

*with inverse given by*

$$(A^*A)^{-1}\alpha = \left[R\alpha(t)\frac{\overline{G}(t)}{M_0(t)} - \int_0^t R\alpha\frac{\overline{G}(u)}{M_0(u)}\frac{\mathrm{d}G}{\overline{G}(u)}\right]g^{1/2}. \tag{37}$$

**Proof.** Using the decompositions $A = A_1 + A_2$ and $A^* = A_1^* + A_2^*$ given in Lemma A.1, write

$$A\alpha = \mu_1 + \mu_2, \qquad \text{where } \mu_i = A_i\alpha, i = 1, 2,$$

so that

$$A^*A\alpha = (A_1^* + A_2^*)(A_1 + A_2)\alpha = A_1^*\mu_1 + A_1^*\mu_2 + A_2^*\mu_1 + A_2^*\mu_2.$$

Consider $A_1^*\mu_1$. From (31) and (32),

$$A_1^*\mu_1 = g^{-1/2}(t)\sum_{\eta,i,r}\int_{z_r} f_0 g^{-1/2}\alpha\,\mathrm{d}z_r$$

$$= g^{-1}(t)\alpha(t)\sum_{\eta,i,r}\int_{z_r} f_0\,\mathrm{d}z_r$$

$$= g^{-1}(t)\alpha(t)\sum_{\eta,i,r}\varrho(\eta)\left[\int_{z_r} K_{\eta,m}g(t)\overline{G}^{\eta-1}(t)h(z_r)\,\mathrm{d}z_r\right]$$

$$= \alpha(t)\sum_{\eta}\varrho(\eta)\left[\overline{G}^{\eta-1}(t)K_{\eta,m}\sum_{i,r}\int_{z_r} h(z_r)\,\mathrm{d}z_r\right]$$

$$= \alpha(t)\sum_{\eta}\varrho(\eta)\left[\overline{G}^{\eta-1}(t)K_{\eta,m}\sum_{i,r} 1\right]$$

$$= \alpha(t)E[\eta\overline{G}^{\eta-1}(t)].$$

In a similar fashion,

$$A_1^*\mu_2 = \sum_{\eta}(\eta-1)g^{-1/2}(t)\sum_{i,r}\int_{z_r} f_0\frac{\int_t^\infty g^{1/2}\alpha}{\overline{G}(t)}\,\mathrm{d}z_r$$

$$= \sum_{\eta}\varrho(\eta)\left[(\eta-1)g^{-1/2}(t)K_{\eta,m}\sum_{i,r}\int_{z_r} g(t)\overline{G}^{\eta-2}(t)\int_t^\infty g^{1/2}\alpha\, h(z_r)\,\mathrm{d}z_r\right]$$

$$= E\left[(\eta-1)g^{1/2}(t)\overline{G}^{\eta-2}(t)\int_t^\infty g^{1/2}\alpha\, K_{\eta,m}\sum_{i,r}\int_{z_r} h(z_r)\,\mathrm{d}z_r\right]$$

$$= E\left[\eta(\eta-1)g^{1/2}(t)\overline{G}^{\eta-2}(t)\int_t^\infty g^{1/2}\alpha\right] = E\left[-\eta(\eta-1)g^{1/2}(t)\overline{G}^{\eta-2}(t)\int_0^t g^{1/2}\alpha\right],$$

where the final equality follows from (3).

Now, moving on to the terms involving $A_2^*$,

$$A_2^*\mu_1 = \sum_\eta (\eta-1)g^{1/2}(t) \sum_{i,r} \int_0^t \int_{z_r} \frac{f_0}{\overline{G}} g^{-1/2}\alpha \, dz_r \, du$$

$$= \sum_\eta \varrho(\eta)\left[(\eta-1)g^{1/2}(t)\int_0^t g^{1/2}\overline{G}^{\eta-2}\alpha \, du K_{\eta,m} \sum_{i,r} \int_{z_r} h(z_r) \, dz_r\right]$$

$$= E\left[\eta(\eta-1)g^{1/2}(t) \int_0^t g^{1/2}\overline{G}^{\eta-2}\alpha \, du\right]$$

and last,

$$A_2^*\mu_2 = \sum_\eta (\eta-1)^2 g^{1/2}(t) \sum_{i,r} \int_0^t \int_{z_r} \frac{f_0}{\overline{G}(u)^2} \int_u^\infty g^{1/2}\alpha \, dz_r \, du$$

$$= \sum_\eta \varrho(\eta)(\eta-1)^2 g^{1/2}(t) \int_0^t g(u)\overline{G}^{\eta-3}(u) \int_u^\infty g^{1/2}\alpha \, K_{\eta,m} \sum_{i,r} \int_{z_r} h(z_r) \, dz_r$$

$$= E\left[-\eta(\eta-1)^2 g^{1/2}(t) \int_0^t g(u)\overline{G}^{\eta-3}(u) \int_0^u g^{1/2}(v)\alpha(v) \, dv \, du\right]$$

$$= E\left[-\eta(\eta-1)^2 g^{1/2}(t) \int_0^t g^{1/2}(v)\alpha(v) \int_v^t g(u)\overline{G}^{\eta-3}(u) \, du \, dv\right]$$

$$= E\left[\frac{-\eta(\eta-1)^2}{\eta-2} g^{1/2}(t) \int_0^t g^{1/2}(v)\alpha(v)[\overline{G}^{\eta-2}(v) - \overline{G}^{\eta-2}(t)]\right]$$

$$= E\left[\frac{\eta(\eta-1)^2}{\eta-2} g^{1/2}(t)\left(\overline{G}^{\eta-2}(t) \int_0^t g^{1/2}\alpha - \int_0^t g^{1/2}\overline{G}^{\eta-2}\alpha\right)\right].$$

Combining terms, we arrive at

$$A^*\mu = E\left[\eta g^{1/2}\left(g^{-1/2}\overline{G}^{\eta-1}\alpha - (\eta-1)\overline{G}^{\eta-2}\int_0^t g^{1/2}\alpha + (\eta-1)\int_0^t g^{1/2}\overline{G}^{\eta-2}\alpha\right.\right.$$
$$\left.\left. + \frac{(\eta-1)^2}{\eta-2}\left(\overline{G}^{\eta-2}\int_0^t g^{1/2}\alpha - \int_0^t g^{1/2}\overline{G}^{\eta-2}\alpha\right)\right)\right]$$

$$= E\left[\eta g^{1/2}\left(g^{-1/2}\overline{G}^{\eta-1}\alpha - (\eta-1)\left(1-\frac{\eta-1}{\eta-2}\right)\left(\overline{G}^{\eta-2}\int_0^t g^{1/2}\alpha - \int_0^t g^{1/2}\overline{G}^{\eta-2}\alpha\right)\right)\right]$$

$$= g^{1/2}E\left[g^{-1/2}\alpha\eta\overline{G}^{\eta-1} + \left(\frac{\eta-1}{\eta-2}\right)\left(\eta\overline{G}^{\eta-2}\int_0^t g^{1/2}\alpha - \int_0^t g^{1/2}\alpha\eta\overline{G}^{\eta-2}\right)\right].$$

Recalling $R\alpha$ from (34), we may write

$$A^*\mu = g^{1/2}E\bigg[g^{-1/2}\alpha\frac{\eta\overline{G}^\eta}{\overline{G}} + \bigg(\frac{\eta-1}{\eta-2}\bigg)\bigg(\frac{\eta\overline{G}^\eta}{\overline{G}^2}\int_0^t g^{1/2}\alpha - \int_0^t g^{-1/2}\alpha\frac{\eta\overline{G}^\eta}{\overline{G}}\frac{\mathrm{d}\overline{G}}{\overline{G}}\bigg)\bigg]$$

$$= g^{1/2}E\bigg[R\alpha\frac{\eta\overline{G}^\eta}{\overline{G}} + \frac{1}{\eta-2}\frac{\eta\overline{G}^\eta}{\overline{G}^2}\int_0^t g^{1/2}\alpha$$

$$- \int_0^t g^{-1/2}\alpha\frac{\eta\overline{G}^\eta}{\overline{G}}\frac{\mathrm{d}\overline{G}}{\overline{G}} - \frac{1}{\eta-2}\int_0^t g^{-1/2}\alpha\frac{\eta\overline{G}^\eta}{\overline{G}}\frac{\mathrm{d}\overline{G}}{\overline{G}}\bigg].$$

Rewriting the third term using

$$\int_0^t g^{-1/2}(u)\alpha(u)\frac{\eta\overline{G}^\eta(u)}{\overline{G}(u)}\frac{\mathrm{d}\overline{G}}{\overline{G}} = \int_0^t \bigg(R\alpha(u) + \frac{\int_u^\infty g^{1/2}\alpha}{\overline{G}}\bigg)\frac{\eta\overline{G}^\eta(u)}{\overline{G}(u)}\frac{\mathrm{d}\overline{G}}{\overline{G}},$$

we find

$$A^*\mu = g^{1/2}E\bigg[R\alpha(t)\frac{\eta\overline{G}^\eta(t)}{\overline{G}(t)} - \int_0^t R\alpha(u)\frac{\eta\overline{G}^\eta(u)}{\overline{G}(u)}\frac{\mathrm{d}\overline{G}}{\overline{G}}$$

$$+ \frac{1}{\eta-2}\bigg(\frac{\eta\overline{G}^\eta}{\overline{G}^2}\int_0^t g^{1/2}\alpha - \int_0^t g^{-1/2}\alpha\frac{\eta\overline{G}^\eta}{\overline{G}}\frac{\mathrm{d}\overline{G}}{\overline{G}} \qquad (38)$$

$$- (\eta-2)\int_0^t \bigg(\int_u^\infty g^{1/2}\alpha\bigg)\frac{\eta\overline{G}^\eta}{\overline{G}^2}\frac{\mathrm{d}\overline{G}}{\overline{G}}\bigg)\bigg].$$

But now we see that the term on second line of (38) vanishes, since

$$(\eta-2)\int_0^t \bigg(\int_u^\infty g^{1/2}\alpha\bigg)\frac{\eta\overline{G}^\eta}{\overline{G}^2}\frac{\mathrm{d}\overline{G}}{\overline{G}} = -\eta\int_0^t\bigg(\int_u^\infty g^{1/2}\alpha\bigg)\mathrm{d}\overline{G}^{\eta-2}$$

$$= -\eta\bigg(\int_u^\infty g^{1/2}\alpha\bigg)\overline{G}^{\eta-2}\bigg|_0^t + \eta\int_0^t \overline{G}^{\eta-2}g^{1/2}\alpha$$

$$= -\eta\overline{G}^{\eta-2}\int_t^\infty g^{1/2}\alpha + \int_0^t g^{1/2}\alpha\eta\overline{G}^{\eta-2}$$

$$= \eta\overline{G}^{\eta-2}\int_0^t g^{1/2}\alpha - \int_0^t g^{-1/2}\alpha\eta\overline{G}^{\eta-2}\,\mathrm{d}\overline{G}$$

$$= \frac{\eta\overline{G}^\eta}{\overline{G}^2}\int_0^t g^{1/2}\alpha - \int_0^t g^{-1/2}\alpha\frac{\eta\overline{G}^\eta}{\overline{G}}\frac{\mathrm{d}\overline{G}}{\overline{G}}.$$

Hence, $A^*A\alpha$ is given by the first line of (38), and taking the expectation inside the integral completes the proof of (36).

Finally, as $A^*A$ is of the form (36), the form (37) of the inverse follows as in [1], page 449. $\square$

We are now in position to prove Lemma 3.2, giving the solution $\hat{\alpha}$ to the normal equations $A^*A\alpha = A^*\rho_0$, and the projection of $\rho_0$ onto the range of $A$.

**Proof of Lemma 3.2.** With $\rho_0$ as in (23), we first claim

$$A^*\rho_0 = \frac{EZ}{2}g^{1/2}(t)E\left[\frac{\eta}{\eta-1}(\eta\overline{G}(t)^{\eta-1} - 1)\right]. \tag{39}$$

From (32) we obtain directly that

$$A_1^*\rho_0 = \frac{EZ}{2}g^{1/2}(t)E[\eta(1+\eta\log\overline{G}(t))\overline{G}^{\eta-1}(t)]$$

and

$$A_2^*\rho_0 = \frac{EZ}{2}g^{1/2}(t)E\left[\eta\int_0^t ((1+\eta\log\overline{G}(u))(\eta-1)g(u)\overline{G}^{\eta-2}(u))\,\mathrm{d}u\right]$$

$$= \frac{EZ}{2}g^{1/2}(t)E\left[\eta\left(\frac{1}{\eta-1}\overline{G}^{\eta-1}(u) - \eta\log\overline{G}(u)\overline{G}^{\eta-1}(u)\right)_0^t\right]$$

$$= \frac{EZ}{2}g^{1/2}(t)E\left[\eta\left(-\frac{1}{\eta-1} + \frac{1}{\eta-1}\overline{G}^{\eta-1}(t) - \eta\log\overline{G}(t)\overline{G}^{\eta-1}(t)\right)\right]$$

and adding these two contributions yields the result (39).

From (34) and (39)

$$R(A^*\rho_0) = \frac{EZ}{2}E\left[\frac{\eta}{\eta-1}\left(\eta\overline{G}(t)^{\eta-1} - 1 - \frac{1}{\overline{G}(t)}\int_t^\infty [\eta\overline{G}^{\eta-1} - 1]\,\mathrm{d}G\right)\right]$$

$$= \frac{EZ}{2}E\left[\frac{\eta}{\eta-1}\left(\eta\overline{G}(t)^{\eta-1} - 1 + \frac{1}{\overline{G}(t)}(-\overline{G}^\eta(t) + \overline{G}(t))\right)\right]$$

$$= \frac{EZ}{2}E[\eta\overline{G}(t)^{\eta-1}] = \frac{1}{2}\frac{M_1(t)}{\overline{G}(t)},$$

where $M_1(t) = E[Z]E[\eta\overline{G}(t)^\eta]$ in accordance with (35).

Hence, by (37), the solution $\hat{\alpha}$ to the normal equations $A^*A\alpha = A^*\rho_0$ is given by

$$\hat{\alpha} = (A^*A)^{-1}A^*\rho_0(t) = \frac{1}{2}\left[\frac{M_1(t)}{M_0(t)} - \int_0^t \frac{M_1(s)}{M_0(s)}\frac{\mathrm{d}G}{\overline{G}}\right]g^{1/2}(t)$$

$$= \frac{1}{2}\left[E(Z) - \int_0^t E(Z)\frac{\mathrm{d}G}{\overline{G}(s)}\right]g^{1/2}(t)$$

$$= \frac{E(Z)}{2}[1 + \log\overline{G}(t)]g^{1/2}(t),$$

where we have used $M_1(t)/M_0(t) = EZ$. To calculate the projection $A\hat{\alpha}$ of $\rho_0$ onto the range of $A$, note

$$\int_t^\infty \hat{\alpha}(s)g^{1/2}(s)\,\mathrm{d}s = \frac{E(Z)}{2}\int_t^\infty (1+\log\overline{G}(s))\,\mathrm{d}G(s) = \frac{E(Z)}{2}\overline{G}(t)\log\overline{G}(t),$$

and hence

$$A\hat{\alpha} = \left[g^{-1/2}(t)\hat{\alpha}(t) + (\eta-1)\frac{\int_t^\infty g^{1/2}\hat{\alpha}}{\overline{G}(t)}\right]f_0^{1/2} = \frac{E(Z)}{2}[1+\eta\log\overline{G}(t)]f_0^{1/2}.$$

$\square$

## A.2. $B$ operator

In this section we prove Lemma 3.4, providing the solution to the normal equations for the operator $B$. Parallel to Section A.1, we begin by deriving an expression for $B^*B$.

**Lemma A.3.** *Let the operator $B$ be given by (25). Then*

$$B^*B\beta = m\beta(z).$$

**Proof.** Applying formulas (27), (25) and (18),

$$B^*B\beta = h^{-1/2}(z)\sum_{\eta,i,r,j\in r}\int_{z_{\neg j}}\int_0^\infty f_0\left(\sum_{k\in r}h^{-1/2}(z_k)\beta(z_k)\right)\mathrm{d}t\,\mathrm{d}z_{\neg j}$$

$$= \sum_\eta \varrho(\eta)K_{\eta,m}h^{-1/2}(z)$$

$$\times \sum_{i,r,j\in r}\int_0^\infty g(t)\overline{G}(t)^{\eta-1}\,\mathrm{d}t\int_{z_{\neg j}}\prod_{l\in r}h(z_l)\left(\sum_{k\in r}h^{-1/2}(z_k)\beta(z_k)\right)\mathrm{d}z_{\neg j}$$

$$= h^{-1/2}(z)E\left[\int_0^\infty \eta g(t)\overline{G}(t)^{\eta-1}\,\mathrm{d}t\right]\int_{z_{\neg j}}\sum_{j\in[m]}\prod_{l\in[m]}h(z_l)\left(\sum_{k\in[m]}h^{-1/2}(z_k)\beta(z_k)\right)\mathrm{d}z_{\neg j}$$

$$= h^{-1/2}(z)\sum_{j\in[m]}\int_{z_{\neg j}}\prod_{l\in[m]}h(z_l)\left(\sum_{k\in[m]}h^{-1/2}(z_k)\beta(z_k)\right)\mathrm{d}z_{\neg j}$$

$$= mh^{-1/2}(z)\int_{z_{\neg 1}}\prod_{l\in[m]}h(z_l)\left(\sum_{k\in[m]}h^{-1/2}(z_k)\beta(z_k)\right)\mathrm{d}z_{\neg 1},$$

where the third equality is by symmetry, and the last by recalling that $z_1$ and $z$ are identified in the integral over $z_{\neg 1}$. Hence

$$B^*B\beta = mh^{1/2}(z)\int_{z_{\neg 1}} \prod_{l=2}^m h(z_l)\left(\sum_{k=1}^m h^{-1/2}(z_k)\beta(z_k)\right) dz_{\neg 1}$$

$$= mh^{1/2}(z)\int_{z_{\neg 1}} \prod_{l=2}^m h(z_l)\left(h^{-1/2}(z)\beta(z) + \sum_{k=2}^m h^{-1/2}(z_k)\beta(z_k)\right) dz_{\neg 1}$$

$$= m\beta(z)\int_{z_{\neg 1}} \prod_{l=2}^m h(z_l)\, dz_{\neg 1} + mh^{1/2}(z)\int_{z_{\neg 1}} \prod_{l=2}^m h(z_l)\left(\sum_{k=2}^m h^{-1/2}(z_k)\beta(z_k)\right) dz_{\neg 1}.$$

As $h(z_l)$ is a density, the first term integrates to $m\beta(z)$. For the second term,

$$\int_{z_{\neg 1}} \prod_{l=2}^m h(z_l) \sum_{k=2}^m h^{-1/2}(z_k)\beta(z_k)\, dz_{\neg 1}$$

$$= \sum_{k=2}^m \int_{z_{\neg 1}} \left(\prod_{l\notin\{1,k\}} h(z_l)\right) h^{1/2}(z_k)\beta(z_k)\, dz_{\neg 1}$$

$$= \sum_{k=2}^m \int_{z_{\neg 1,k}} \left(\prod_{l\notin\{1,k\}} h(z_l)\right) dz_{\neg 1,k} \int_{z_k} h^{1/2}(z_k)\beta(z_k)\, dz_k$$

$$= \sum_{k=2}^m \int_{z_k} h^{1/2}(z_k)\beta(z_k)\, dz_k$$

$$= 0$$

by (6), showing $B^*B\beta = m\beta(z)$ and the lemma. $\square$

**Proof of Lemma 3.4.** From (27) and (23), arguing as in the proof of Lemma A.3 and applying Lemma 4.3, we obtain

$$B^*\rho_0 = \frac{1}{2}\sum_{\eta,i,r,j\in r} \int_{z_{\neg j}}\int_0^\infty f_0 h^{-1/2}(z_j)\left(z_i + \log\overline{G}(t)\sum_{k\in r}(z_k - EZ) + \eta\log\overline{G}(t)EZ\right)dt\, dz_{\neg j}$$

$$= \frac{1}{2}\sum_\eta \varrho(\eta)\left[\sum_{j\in[m]} \int_{z_{\neg j}}\int_0^\infty \eta g(t)\overline{G}^{\eta-1}(t)h^{-1/2}(z_j)h(z_r)\right.$$

$$\left.\times \left(z_1 + \log\overline{G}(t)\sum_{k\in[m]}(z_k - EZ) + \eta\log\overline{G}(t)EZ\right)dt\, dz_{\neg j}\right]$$

$$= \frac{1}{2} \sum_{j \in [m]} \int_{z_{\neg j}} h^{-1/2}(z_j) h(z_r) E\left[z_1 - \frac{1}{\eta} \sum_{k \in [m]} (z_k - EZ) - EZ\right] \mathrm{d}z_{\neg j}$$

$$= \frac{1}{2} \sum_{j \in [m]} \int_{z_{\neg j}} h^{1/2}(z_j) \prod_{k \neq j} h(z_k) E\left[\left(\frac{\eta-1}{\eta}\right)(z_1 - EZ) - \frac{1}{\eta} \sum_{k=2}^{m} (z_k - EZ)\right] \mathrm{d}z_{\neg j}.$$

For $j = 1$ we obtain $(1/2)h^{1/2}(z)E[(\eta-1)/\eta](z-EZ)$ from the first term in parentheses, while each term in the second sum integrates to zero. For each of the $m-1$ terms, where $j \neq 1$ the first term in parentheses integrates to zero, but when $k = j$ one term in the sum in the second term makes a non-zero contribution of $-(1/2)h^{1/2}(z)(z-EZ)E[1/\eta]$, for a total of

$$B^* \rho_0 = \frac{1}{2} h^{1/2}(z)(z-EZ)E\left[\frac{\eta-1}{\eta} - \frac{m-1}{\eta}\right] = \frac{1}{2} h^{1/2}(z)(z-EZ)E\left[\frac{\eta-m}{\eta}\right].$$

From Lemma A.3 we clearly have

$$(B^*B)^{-1}\beta = \frac{1}{m}\beta \quad \text{hence} \quad \hat{\beta} = (B^*B)^{-1}B^*\rho_0 = \frac{1}{2}h^{1/2}(z)(z-EZ)E\left[\frac{\eta-m}{\eta m}\right],$$

proving (28). Applying $B$ as in (25) to $\hat{\beta}$ now yields (29). □

## A.3. Verification of Assumptions 2.1–2.3

In this section we provide a basic outline of the verifications of Assumptions 2.1–2.3 given in detail in the technical report [9]. In particular, it is shown there by a simple argument that Assumption 2.1 is satisfied with

$$\mathcal{A} = \{\alpha \in L^2(\nu^+) : \langle \alpha, g^{1/2} \rangle_{\nu^+} = 0\} \quad \text{and} \quad \mathcal{B} = \bigcap_{|\theta| < \theta_\kappa} \{\beta \in L^2(\nu_\theta) : \langle \beta, h^{1/2} \rangle_\nu = 0\}. \quad (40)$$

Given the quantities $A, B, \hat{\alpha}, \hat{\beta}$ and $\rho_0$ in (24), (25), Lemmas 3.2, 3.4 and (23), respectively, it is lengthy to verify the remainder of Assumption 2.2, that is, the required convergence (8). One main point of the detailed verification given in [9] is that, with

$$\overline{G}_{n,\theta}(t;z) = \overline{G}_n^{e^{\theta z}}(t), \qquad \overline{G}_{n,\theta}(t) = \int \overline{G}_{n,\theta}(t;z) h_n(z) \, \mathrm{d}z \quad \text{and} \quad \overline{G}_n(t) = \overline{G}_{n,0}(t),$$

each of the three conditions given in Theorem 2.5 implies that

$$C_n = \left(\frac{\eta-m}{2\eta}\right) \frac{E_n[Z e^{\theta Z} \overline{G}_n(t)^{e^{\theta z}}]}{\overline{G}_{n,\theta}(t)}$$

is uniformly bounded.

The positivity condition is the one that appears in [3]. Under positivity $z \geq 0$ and $\theta \geq 0$, the function $z e^{\theta z}$ is increasing in $z$ and $\overline{G}_n(t)^{e^{\theta z}}$ is decreasing in $z$. Therefore these functions are negatively correlated and we have

$$E_n[Z e^{\theta Z} \overline{G}_n(t)^{e^{\theta z}}] \leq E_n[Z e^{\theta Z}] E_n[\overline{G}_n(t)^{e^{\theta z}}] = E_n[Z e^{\theta Z}] \overline{G}_{n,\theta}(t),$$

and in particular

$$|C_n| \leq \tfrac{1}{2} E_n[Z e^{\theta Z}],$$

which is a bounded sequence in $n$. Under the bounded covariate condition with, say, $|Z| \leq z_0$ almost surely, we have

$$|C_n| \leq \frac{1}{2} \frac{E_n[|Z e^{\theta Z}| \overline{G}_n(t)^{e^{\theta z}}]}{E_n[\overline{G}_n(t)^{e^{\theta z}}]} \leq \frac{1}{2} z_0 e^{|\theta| z_0} \leq \frac{1}{2} z_0 e^{\overline{\theta} z_0}.$$

Under the cohort size condition one shows that

$$C_n = \left(\frac{\eta - m}{2\eta}\right) E_n[Z e^{\theta Z} \overline{G}_n(t)^{e^{\theta Z}}] \left(\frac{K_{\eta,m}}{K_{\eta-2,m}}\right)^{1/2},$$

and since $|E_n[Z e^{\theta Z} \overline{G}_n(t)^{e^{\theta Z}}]| \leq E_n[|Z| e^{\theta Z}]$ is bounded in $n$, and $K_{\eta,m}/K_{\eta-2,m} \leq 1$, again we find the constant $C_n$ to be uniformly bounded.

Regarding Assumption 2.3, by (40) to show that, say, $\hat{\alpha} \in \mathcal{A}$, it suffices to verify that $\hat{\alpha} \in L^2(\nu^+)$ and $\langle \hat{\alpha}, g^{1/2} \rangle_{\nu^+} = 0$. The first claim follows from Lemma 4.3, and, by applying that same lemma with $\eta = k$ and $j = 1$, the second claim from

$$\int_0^\infty \log \overline{G}(t) g(t) \, dt = -1.$$

The verification that $\hat{\beta} \in \mathcal{B}$ is similar, but somewhat more involved.

### A.4. Proofs of Theorems 2.1–2.4

**Proof of Theorem 2.1.** The set $\mathbb{H}$ given in (10) is a subspace of $L^2(\sigma)$, being the image of the subspace $\mathbb{R} \times \mathcal{A} \times \mathcal{B}$ under the linear transformation $(\tau, \alpha, \beta) \to \tau \rho_\theta + A\alpha + B\beta$. Hence, with $\hat{\alpha}$ and $\hat{\beta}$ as in Assumption 2.2, by that assumption,

$$\tau \hat{\zeta} \in \mathbb{H} \qquad \text{for all } \tau \in \mathbb{R}, \text{ where } \hat{\zeta} = \rho_\theta - A\hat{\alpha} - B\hat{\beta}, \tag{41}$$

and $4\|\tau\hat{\zeta}\|_\sigma^2 = \tau^2 I_*$. Now let $\{f_n\}_{n \geq 0} \in \mathcal{F}(f, \tau\hat{\zeta})$ and continue as in the proof of Theorem 3.1 in [1]. In particular, with $L_n$ as the log likelihood ratio for $f_n$ vs. $f_0$, $S$ as the limiting distribution of $n^{1/2}(\hat{\theta}_n - \theta_0)$, guaranteed to exist by the regularity of $\hat{\theta}_n$, and $Z \sim \mathcal{N}(0, I_*)$, the random vector $(n^{1/2}(\hat{\theta}_n - \theta), L_n)$ converges weakly under $f$ to $(S, \tau Z - 1/2\tau^2 I_*)$ and the characteristic function of $S$ factors into the product of the characteristic functions of $S - Z/I_*$ times that of $Z/I_*$. $\square$

**Proof of Theorem 2.2.** Note that with $\hat{\zeta}$ as in (41), as $\hat{\zeta} \in \mathbb{H}$ we have that $\mathcal{F}(f,\hat{\zeta}) \subset \mathcal{F}(f)$, and therefore, for all $c \geq 0$

$$B_n^*(c) = \{f_n \in \mathcal{F}(f,\hat{\zeta}) : n^{1/2}\|f_n^{1/2} - f^{1/2}\|_\sigma \leq c\}$$
$$\subset \{f_n \in \mathcal{F}(f) : n^{1/2}\|f_n^{1/2} - f^{1/2}\|_\sigma \leq c\} = B_n(c).$$

Hence the argument for the proof of Theorem 3.2 of [1] is obtained. $\square$

**Proofs of Theorems 2.3 and 2.4.** The application of the results of [13] and [2], as in the proofs of Theorems 4.1 and 4.2 in [1] apply with minimal changes. In particular, for any element of $\mathbb{H}$ given by

$$\zeta = \tau\rho_\theta + A\alpha + B\beta,$$

let $\mathcal{T} : H \to B_0 = \{x \in C[0,1] : x(0) = x(1) = 0\}$ be defined by

$$\mathcal{T}\zeta = \int_0^t 2\alpha g^{1/2}\,\mathrm{d}\nu^+ = \langle \alpha, 2g^{1/2}I_{[0,t]}\rangle_{\nu^+}.$$

Writing $\zeta$ as

$$\zeta = \tau(\rho_\theta - A\hat{\alpha} - B\hat{\beta}) + A(\tau\hat{\alpha} + \alpha) + B(\tau\hat{\beta} + \beta),$$

the orthogonality provided by Assumption 2.3 and Lemma 4.1 yield $A^*\zeta = A^*A(\tau\hat{\alpha} + \alpha)$ and therefore

$$\alpha = C^*\zeta - \langle\zeta, 4(\rho_\theta - A\hat{\alpha} - B\hat{\beta})/I_*\rangle\hat{\alpha}, \qquad \text{where } C = A(A^*A)^{-1}.$$

Continuing, one may verify that the adjoint $\mathcal{T}^*$ of $\mathcal{T}$ is given by a formula analogous to that in Lemma 5.2 of [1], and that

$$\frac{1}{4}\|\mathcal{T}^*v\|_\sigma^2 = E\left(\int_0^t \mathbb{Z}_*\,\mathrm{d}\nu\right)^2.$$

$\square$

We remark that though the subspace $\mathbb{H}$ is not assumed to be closed in $L^2(\sigma)$, and hence the projection theorem cannot be applied, as long as $\mathbb{H}$ contains the approach to $f$ along the 'worst case' direction $\hat{\zeta}$, the proof of [1] carries through. Moreover, this holds true independently of the number of factors in the model, one more here than in [1]. The other difference between the situation here and that of [1], that $\mathcal{B}$ consists of the perturbations that approximate $n^{1/2}(h_n^{1/2} - h^{1/2})$ in $L^2(\nu_\theta)$ for all $|\theta| < \theta_\kappa$ rather than in the weaker $L^2(\nu)$ sense, is handled by Assumption 2.2, which gives, in particular, that the critical $\hat{\beta}$ lies in $\mathcal{B}$ even when insisting on the stronger form of convergence.

598                                                                    *L. Goldstein and H. Zhang*

# Acknowledgement

The authors acknowledge the support of National Cancer Institute Grant CA 42949.